%% file: Otnositel_naya_dissipatsiya.tex
\documentclass[11pt,reqno]{amsart}
\usepackage{amssymb,amscd,amsbsy}
\usepackage{amssymb,amscd,amsbsy,mathrsfs}
\input{classstartscr}

\theoremstyle{remark}

\newtheorem*{rem*}{Remark}

\newcommand{\OLA}{{\rm OL}_{\rm A}}
\newcommand{\OL}{{\rm OL}}
\newcommand{\ROL}{{\rm ROL}}
\newcommand{\OLAr}{{\rm OL}_{\rm A}^{\rm res}}
\newcommand{\COL}{{\rm CL}}

\newcommand{\MOL}{{\frak M}_{\OL}}
\newcommand{\MOLA}{{\frak M}_{\OL_A}}
\newcommand{\CAb}{{\rm C}_{{\rm A},\be}}
\newcommand{\CAr}{{\rm C}_{\rm A}(\C_+)}
\newcommand\mM{\mathcal{M}}

\newcommand{\Dom}{{\rm Dom}}

\newcommand\fM{\frak M}

\newcommand\dg{\frak D}

\newcommand{\rd}{{\rm d}}

\newcommand\ri{{\rm i}}

\begin{document}

\newcommand{\vse}{\vspace{.2in}}
\numberwithin{equation}{section}

\author{A.B. Aleksandrov and V.V. Peller}

\title{Functions of dissipative operators under relatively bounded and relatively trace class perturbations}
\thanks{The research on \S\:4-6  is supported by 
Russian Science Foundation [grant number 23-11-00153].
The research on \S\:7-10 is supported by a grant of the Government of the Russian Federation for the state support of scientific research, carried out under the supervision of leading scientists, agreement  075-15-2025-013.}

\begin{abstract}
We study the behaviour of functions of dissipative operators under relatively bounded and relatively trace class perturbation. We introduce and study the class of analytic relatively operator Lipschitz functions. An essential role is played by double operator integrals with respect to semispectral measures. We also study the class of analytic resolvent Lipschitz functions. 
Then we obtain a trace formula in the case of relatively trace class perturbations and show that the maximal class of function for which the trace formula holds in the case of relatively trace class perturbations coincides with the class of analytic relatively operator Lipschitz functions. We also establish the inequality 
$\int|\boldsymbol{\xi}(t)|(1+|t|)^{-1}\,{\rm d}t<\be$ for the spectral shift function $\boldsymbol\xi$  in the case of relatively trace class perturbations.
\end{abstract}

\maketitle

{\bf
\footnotesize
\tableofcontents
\normalsize
}

\setcounter{section}{0}
\section{\bf Introduction}
\setcounter{equation}{0}
\label{In}

\

In our  recent paper \cite{AP6} (see also \cite{AP7}) we studied the behaviour of functions of self-adjoint operators under relatively bounded and relatively compact perturbations. Recall that a linear operator $K$ is called a {\it relatively bounded perturbation} of an unbounded self-adjoint operator $A$ if the operator 
$K(A+\ri I)^{-1}$ is bounded and it is called a {\it relatively trace class perturbation} of $A$ if the operator $K(A+\ri I)^{-1}$ belongs to the trace class $\bS_1$. We introduced in \cite{AP6} and \cite{AP7} the class $\ROL$
of {\it relatively operator Lipschitz functions} as the class of functions $f$ on the real line $\R$, for which tha inequality
$$
\|f(B)-f(A)\|\le k\|(B-A)(A+\ri I)^{-1}\|
$$
holds for arbitrary self-adjoint operators $A$ and $B$ such that $B-A$ is a relatively bounded perturbation of $A$.

We found in \cite{AP6} and \cite{AP7} various characterizations of the class $\ROL$ of relatively operator Lipschitz functions. 

We also obtained in \cite{AP6} and \cite{AP7} the trace formula
\bay
\label{sledy}
\trace\big(f(B)-f(A)\big)=\int_\R f'(t)\bs\xi(t)\,\rd t
\ey
for arbitrary $f\in\ROL$ and arbitrary pairs if self-adjoint operators $A$ and $B$ such that $B-A$ is a relatively trace class perturbation of $A$. Here $\bs\xi$ is the {\it spectral shift function} for the pair $\{A,B\}$, which satisfies the condition
$$
\int\frac{|\bs\xi(t)|}{1+|t|}<\be
$$
(the last inequality was also obtained in \cite{CS}). Moreover, it turned out that $\ROL$ is the maximal class of functions, for which the trace formula \rf{sledy} holds for arbitrary pairs $\{A,B\}$ such that $B-A$ is a relatively trace class perturbation of $A$, see \cite{AP6} and \cite{AP7}.

In this paper we consider similar problems for functions of dissipative operators. Note, however, that certain proofs in the dissipative case are considerably harder than in the self-adjoint case.

\

\section{\bf Maximal dissipative operators. Functional calculus}
\setcounter{equation}{0}
\label{maks}

\

Let $\h$ be a Hilbert space. An operator $L$ (not necessarily bounded)
with dense domain $\Dom(L)$ in $\h$ is called {\it dissipative} if
$$
\im(Lu,u)\ge0,\quad u\in\Dom(L).
$$
A dissipative operator is called {\it maximal dissipative} if it has no proper dissipative extension.


The {\it Cayley transform} of a dissipative operator $L$ is defined by
$$
T\df(L-{\rm i}I)(L+{\rm i}I)^{-1}
$$
with domain $\Dom(T)=(L+{\rm i}I)\Dom(L)$ and range $\Range T=(L-{\rm i}I)\Dom(L)$
(the operator $T$ is not densely defined in general). $T$ is a contraction, i.e. $\|Tu\|\le\|u\|$, $u\in\Dom(T)$, $1$ is not an eigenvalue of $T$, and $\Range(I-T)\df\{u-Tu:~u\in\Dom(T)\}$ is dense.

Conversely, if $T$ is a contraction defined on its domain $\Dom(T)$ such that $1$ is not an eigenvalue of $T$ and $\Range(I-T)$ is dense, then it is the Cayley transform of a dissipative operator $L$ and $L$ is the inverse Cayley transform of $T$:
$$
L={\rm i}(I+T)(I-T)^{-1},\quad\Dom(L)=\Range(I-T).
$$

A dissipative operator is maximal if and only if the domain of its Cayley transform is the whole Hilbert space.
Clearly, a maximal dissipative operator must be closed.

If $L$ is a dissipative operator, which is not maximal, there is a canonical way to extend $L$ to a maximal dissipative operator. Indeed, let $T$ be the Cayley transform of $L$. Then $T$ extends uniquely to the closure of its domain. If the closure is a proper subspace of our Hilbert space, we can extend $T$ to the whole space by defining it to be the zero operator on the orthogonal complement to its domain. The inverse Cayley transform of this extension is a maximal dissipative operator. It is natural to consider it as the {\it canonical extension of $L$ to a maximal dissipative operator}. 

We proceed now to the construction of functional calculus for dissipative operators. Let $L$ be a maximal dissipative operator and let $T$ be its Cayley transform. Consider its minimal unitary dilation $U$, i.e. $U$ is a unitary operator defined on a Hilbert space $\K$ that contains $\h$ such that
$$
T^n=P_\h U^n\big|\h,\quad n\ge0,
$$
and $\K=\clos\spn\{U^nh:~h\in\h,~n\in\Z\}$. Since $1$ is not an eigenvalue of $T$, it follows that $1$ is not an eigenvalue of $U$ (see \cite{SNF}, Ch. II, \S\,6).

The Sz.-Nagy--Foia\c s functional calculus allows us to define a functional calculus for $T$ on the Banach algebra
$$
{\rm C}_{{\rm A},1}\df\big\{g\in H^\be:~g\quad
\mbox{is continuous on}\quad\T\setminus\{1\}~\big\}.
$$
If $g\in{\rm C}_{{\rm A},1}$, we put
$$
g(T)\df P_\h g(U)\Big|\h.
$$
This functional calculus is linear and multiplicative, and
$$
\|g(T)\|\le\|g\|_{H^\be},\quad g\in{\rm C}_{{\rm A},1}.
$$

We can define now a functional calculus for our dissipative operator on the
Banach algebra 
$$
{\rm C}_{{\rm A},\be}\df\big\{f\in H^\be(\C_+):~f\quad\mbox{is continuous on}\quad
\R\big\}.
$$
Indeed, if $f\in{\rm C}_{{\rm A},\be}$, we put
$$
f(L)\df \big(f\circ\o\big)(T),
$$
where $\o$ is the conformal map of $\dd$ onto $\C_+$ defined by
$\o(\z)\df{\rm i}(1+\z)(1-\z)^{-1}$, $\z\in\dd$.

Let us define now a resolvent self-adjoint dilation of a maximal dissipative operator. If $L$ is a maximal dissipative operator in a Hilbert space $\h$, we say that a self-adjoint operator $A$ in a Hilbert space
$\K$, $\K\supset\h$, is called a {\it resolvent self-adjoint dilation} of $L$
if
$$
(L-\l I)^{-1}=P_\h(A-\l I)^{-1}\Big|\h,\quad \im\l<0.
$$
The dilation is called {\it minimal} if
$$
\K=\clos\spn\big\{(A-\l I)^{-1}v:~v\in\h,~\im\l<0\big\}.
$$
If $f\in\CAb$, then
$$
f(L)=P_\h f(A)\Big|\h,\quad f\in\CAb.
$$

A minimal resolvent self-adjoint dilation of a maximal dissipative operator always exists (and is unique up to a natural isomorphism). Indeed, it suffices to take a minimal unitary dilation of the Cayley transform of this operator and apply the inverse Cayley transform to it. 

We can define now the semi-spectral measure $\E_L$ of a maximal dissipative operator $L$ by
$$
\E_L(\D)=P_\h E_A(\D)\big|\h,
\quad\D\quad\mbox{is a Borel subset of}\quad\R,
$$
where $E_A$ is the spectral measure of the minimal resolvent self-adjoint dilation of $L$.
Then
\bay
\label{fL}
f(L)=\int_\R f(x)\,d\E_L(x),\quad f\in\CAb.
\ey

Note that the semi-spectral measure of a maximal dissipative operator $L$ can be obtained from the semi-spectral measure of the Cayley transform of $L$ by conformally transferring it from the unit circle to the real line.

We refer the reader to \cite{AP1} for more detail.

We also need functions $f(L)$ of maximal dissipative operators in the case when $f$ is not necessarily
a bounded function.
Suppose that $f$ is a  Lipschitz function in $\clos\C_+$ that is analytic in $\C_+$. We can define the operator $f(L)$ as follows. We have
\bay
\label{fri}
f(\z)=\frac{f_\ri(\z)}{(\z+\ri)^{-1}},\quad\z\in\clos\C_+,
\qquad
\mbox{where}\quad f_\ri(\z)\df\frac{f(\z)}{\z+\ri}.
\ey
Clearly, $f_\ri\in{\rm C}_{{\rm A},\be}$
The (possibly unbounded) operator $f(L)$ can be defined by
\bay
\label{f(L)}
f(L)\df(L+\ri I)f_\ri(L)
\ey
(see \cite{SNF}, Ch. IV, \S\;1). It follows from Th. 1.1 of Ch. IV of \cite{SNF}  that 
$$
f(L)\supset f_\ri(L)(L+\ri I),
$$
and so $\Dom(f(L))\supset\Dom(L)$.

We are ready now to define operator Lipschitz functions analytic in the upper half-plane $\clos\C_+$.

\medskip

{\bf Definition.} Suppose that $f$ is a Lipschitz function in the closed upper half-plane $\clos\C_+$ that is analytic in the open half-plane $\C_+$. We say that $f$ is an {\it operator Lipschitz function} if
$$
\|f(L)-f(M)\|\le\const\|L-M\|,
$$
whenever $L$ and $M$ are maximal dissipative operators with bounded difference. 

\medskip

Note that in the case when $L-M$ is bounded, both operators $f(L)$ and $f(M)$ are defined on the dense subset $\Dom(L)=\Dom(M)$.

We denote by 
$\OL_{\rm A}$ the class of such operator Lipschitz functions.

Note also that a function $f$ analytic in $\C_+$ and continuous in $\clos\C_+$ is operator Lipschitz if and only if
$$
\|f(N_1)-f(N_2)\|\le\const\|N_1-N_2\|,
$$
whenever $N_1$ and $N_2$ are bounded normal operators with spectra in $\clos\C_+$, see \cite{APol} and
\cite{AP4}.
We refer the reader to the survey \cite{APol} for detailed information on operator Lipschitz functions.

In \S\:\ref{SlROLA} of this paper we obtain a trace formula in the case when $L$ and $M$ are maximal dissipative operators such that $M-L$ is a relatively trace class perturbation of $L$.

\

\section{\bf Trace formulae}
\setcounter{equation}{0}
\label{Sledy}

\

Let $A$ be a (not necessarily bounded) self-adjoint operator and $B=A+K$, where $K$ is a trace class self-adjoint operator that is considered as a perturbation of $A$. Physicist I.M. Lifshits in \cite{L} proposed the following trace formula
\bay
\label{LifKr}
\trace\big(f(B)-f(A)\big)=\int_\R f'(t)\bs\xi(t)\,\rd t
\ey
for sufficiently nice functions $f$ on $\R$, where $\bs\xi$ is an integrable function on $\R$. It is called the {\it spectral shift function associated with} the pair $\{A,B\}$. Later M.G. Krein found in \cite{Kr} a rigorous mathematical justification of trace formula \rf{LifKr}. We are going to refer to \rf{LifKr} as the {\it Lifshits--Krein trace formula}. Krein posed in \cite{Kr} the problem to identify the maximal class of functions $f$ on $\R$, for which formula \rf{LifKr} holds for all pairs $\{A,B\}$ of self-adjoint operators with trace class difference. Krein's problem was solved in the paper \cite{PeKr}: the maximal class of such functions $f$ coincides with the class $\OL$ of operator Lipschitz functions on $\R$. 

If we replace the condition $B-A\in\bS_1$ with the resolvent condition
\bay
\label{rezus}
(B+\ri I)^{-1}-(A+\ri I)^{-1}\in\bS_1,
\ey
then trace formula \rf{LifKr} holds for the class of {\it resolvent Lipschitz functions} $f$, i.e. for functions $f$
satisfying the inequality
$$
\|f(B)-f(A)\|\le\const\big\|(B+\ri I)^{-1}-(A+\ri I)^{-1}\big\|
$$
for arbitrary pairs $\{A,B\}$ of self-adjoint operators. The spectral shift function $\bs\xi$ under condition
\rf{rezus} must satisfy the inequality
\bay
\label{1+t2}
\int_\R\frac{|\bs\xi(t)|}{1+t^2}\rm\rd t<\be.
\ey
This is a well-known fact and can be reduced to the trace formula for unitary operators, see \cite{MNP} for details. The maximal class of functions $f$, for which trace formula \rf{LifKr} holds for arbitrary pairs $\{A,B\}$ of self-adjoint operators satisfying the resolvent condition 
\rf{rezus} coincides with the class of resolvent Lipschitz functions. The latter follows from the corresponding result for functions of unitary operators \cite{APun}, see \cite{MNP} for details. We refer the reader to \cite{AP6} and \cite{AP7} for various characterizations of the class of resolvent Lipschitz functions.

Finally, consider the situation when $K=B-A$ is a relatively trace class perturbation of $A$. In this case trace formula \rf{LifKr} holds for arbitrary relatively operator Lipschitz functions and the spectral shift function $\bs\xi$ in this situation must satisfy the inequality
$$
\int_\R\frac{|\bs\xi(t)|}{1+|t|}\rm\rd t<\be,
$$
see \cite{CS} and \cite{AP7} (see also \cite{AP6}).
Moreover,
the class $\ROL$ of relatively operator Lipschitz functions is the maximal class of functions, for which trace formula \rf{LifKr} holds for arbitrary pairs 
$\{A,B\}$ of self-adjoint operators such that $B-A$ is a relatively trace class perturbation of $A$. 

We proceed now to trace formulae for functions of dissipative operators. The long-standing problem to generalize the Lifshitz--trace Krein formula to the case of functions of maximal dissipative operators in the case of trace class perturbations and in the case of trace class resolvent perturbations was solved in \cite{MNP1} and \cite{MNP}, se also \cite{MN}. We also refer the reader to \cite{MNP} for the history of 
the question. 

Suppose that $L$ is a maximal dissipative operator and $M$ is a dissipative operator such that 
$M-L\in\bS_1$. Then there exists an integrable function $\bs\xi$ on $\R$ such that
for $f\in\OL_{\rm A}$, the following trace formula holds
\bay
\label{sledis}
\trace\big(f(M)-f(L)\big)=\int_\R f'(t)\bs\xi(t)\,\rd t.
\ey
Such a function $\bs\xi$ (which is never unique) is called a {\it spectral shift function for the pair}
$\{L,M\}$. As in the case of self-adjoint operators, $\OL_{\rm A}$ is the maximal class of functions $f$, for which formula \rf{sledis} holds for all such pairs $\{L,M\}$.

In the case when $\{L,M\}$ is a pair of maximal dissipative operators such that $M$ is a resolvent trace class perturbation of $L$, i.e.
\bay
\label{rezyad}
(M+\ri I)^{-1}-(L+\ri I)^{-1}\in\bS_1,
\ey
trace formula \rf{sledis} holds for a spectral shift function $\bs\xi$ satisfying \rf{1+t2}
and for
functions $f$ in the class $\OLAr$ of {\it analytic resolvent Lipschitz functions}, i.e. for functions $f$ in ${\rm C}_{{\rm A},\be}$ such that
$$
\|f(M)-f(L)\|\le\const\|(M+\ri I)^{-1}-(L+\ri I)^{-1}\|
$$
for arbitrary maximal dissipative operators $L$ and $M$. It was also shown in \cite{MNP} that
$\OLAr$ is the maximal class of functions, for which trace formula \rf{sledis} holds for arbitrary
pairs $\{L,M\}$  of maximal dissipative operators satisfying \rf{rezyad}.
We obtain in \S\:\ref{khara} of this paper various characterizations of the class $\OLAr$.

\

\section{\bf Double operator integrals}
\setcounter{equation}{0}
\label{DOpI}

\

Double operator integrals
\bay
\label{DOpIn}
\iint_{\X\times\Y}\Phi(x,y)\,dE_1(x)Q\,dE_2(y)
\ey
appeared first in the paper \cite{DK} by Yu.L. Daletskii and S.G. Krein. Later
Birman and Solomyak elaborated their beautiful theory of double operator integrals,
see \cite{BS1}, \cite{BS2} and \cite{BS4} (see also \cite{Pe4}, \cite{APol} and references therein). Here $\Phi$ is a bounded measurable function, $E_1$ and $E_2$ are spectral measures on Hilbert space
defined on $\s$-algebras of subsets of sets $\X$ and $\Y$ and
$Q$ is a bounded linear operator. 

The starting point of the Birman--Solomyak approach \cite{BS1} is the case when $Q$ belongs to the Hilbert--Schmidt class $\bS_2$; they managed to define double operator integrals for arbitrary bounded measurable functions $\Phi$. Moreover, for $Q\in\bS_2$ and for a bounded measurable function $\Phi$, the operator
\rf{DOpIn} must belong to $\bS_2$ and the inequality
$$
\left\|\iint_{\X\times\Y}\Phi(x,y)\,dE_1(x)Q\,dE_2(y)\right\|_{\bS_2}
\le\|\Phi\|_{L^\be}\|Q\|_{\bS_2}
$$
must hold.

Next, $\Phi$ is said to be a {\it Schur multiplier with respect to} $E_1$ and $E_2$ if
$$
Q\in\bS_1\quad\Longrightarrow\quad\iint_{\X\times\Y}\Phi(x,y)\,dE_1(x)Q\,dE_2(y)\in\bS_1
$$
(see \cite{Pe1} and \cite{APol}) for more detail). In this case one can define by duality the double operator integral \rf{DOpIn} for an arbitrary bounded linear operator $Q$ and the following inequality must hold
$$
\left\|\iint_{\X\times\Y}\Phi(x,y)\,dE_1(x)Q\,dE_2(y)\right\|
\le\|\Phi\|_{\fM_{E_1,E_2}}\|Q\|,
$$
where $\fM_{E_1,E_2}$ stands for the space of Schur multipliers with respect to $E_1$ and $E_2$ equipped with the natural norm.

It is well known (see \cite{Pe1}, \cite{APol}, \cite{AP5} ) that $\Phi\in\fM_{E_1,E_2}$ if and only if 
$\Phi$ belongs to the {\it Haagerup tensor product} 
$L^\be_{E_1}\otimes_{\rm h}L^\be_{E_2}$ of 
$L^\be_{E_1}$ and $L^\be_{E_2}$ or, in other words, $\Phi$
admits a representation
\bay
\label{htenppre}
\Phi(x,y)=\sum_{n\ge0}\f_n(x)\psi_n(y),
\ey
where the $\f_n$ and $\psi_n$ are measurable functions such that
$$
\sum_{n\ge0}|\f_n|^2\in L^\be_{E_1}\quad\mbox{and}\quad
\sum_{n\ge0}|\psi_n|^2\in L^\be_{E_2}.
$$
In this case
\bay
\label{DOIHTP}
\iint\Phi(x,y)\,dE_1(x)Q\,dE_2(y)=
\sum_{n\ge0}\left(\int\f_n(x)\,dE_1(x)\right)Q\left(\int\psi_n(y)\,dE_2(y)\right),
\ey
the series on the right converges in the weak operator topology
and the right-hand side does not depend on the choice of a representation in 
\rf{htenppre}. Note that the sufficiency of the condition $\Phi\in L^\be_{E_1}\otimes_{\rm h}L^\be_{E_2}$
for $\Phi$ to be a Schur multiplier was established earlier by Birman and Solomyak in \cite{BS1}.

We are going to use the notation $\fM(\R^2)$ the class of functions on $\R^2$ that are Schur multipliers with respect to arbitrary Borel spectral measures on $\R$. 

We also need double operator integrals with respect to {\it semi-spectral measures}
\bay
\label{dvoopipolu}
\iint\Phi(x,y)\,d\E_1(x)Q\,d\E_2(y).
\ey
Such double operator integrals were introduced in \cite{Pe2} (see also \cite{Pe4}).
By analogy with the case of double operator integrals with respect to spectral measures, double operator integrals of the form \rf{dvoopipolu} can be defined for arbitrary bounded measurable functions $\Phi$ in the case when $Q\in\bS_2$ and for 
Schur multipliers $\Phi$ with respect to the semi-spectral measures $\E_1$ and $\E_2$
in the case of an arbitrary bounded operator $Q$. The latter is equivalent to the fact that
$\Phi$ belongs to the Haagerup tensor product
$L^\be_{\E_1}\otimes_{\rm h}L^\be_{\E_2}$ in which case  the definition of \rf{dvoopipolu} is similar to formula \rf{DOIHTP}.

\

\section{\bf Relatively bounded and relatively trace class\\ perturbations of dissipative operators}
\setcounter{equation}{0}
\label{otno}

\

Recall that for a linear operator $A$ in a Hilbert space, its domain is denoted by $\Dom(A)$.
Let $A$ and $K$ be densely defined linear operators such that
$$
\Dom(A)\subset\Dom(K).
$$
It is well known (see \cite{BSSst}, Theorem 1, Ch. 3, Sec. 4) that if the operator $A$ is closed and the operator $K$ is closable, then $K$ is {\it dominated by} $A$, i.e.
\bay
\label{otogner}
\|Kv\|\le\const(\|v\|+\|Av\|),\quad v\in{\rm Dom}(A).
\ey

Suppose now that $L$ and $M$ are maximal dissipative operators such that $K\df M-L$ is dominated by $L$.
In this case we say that $K$ is a {\it relatively bounded perturbation of} $L$. Note that 
$K$ is a relatively bounded perturbation of $L$ if and only if the operator 
\bay
\label{C}
C\df K(L+\ri I)^{-1}
\ey
is bounded. Indeed,
it is obvious that $K(L+\ri I)^{-1}$ is a bounded operator if and only if
$$
\|Kv\|\le\const\|(L+\ri I)v\|,\quad v\in\Dom(L).
$$
Clearly,
$$
\|(L+\ri I)v\|\le\|v\|+\|Lv\|.
$$
It remains to observe that for $v\in\Dom(L)$,
$$
\|(L+\ri I)v\|^2\ge|((L+\ri I)v,(L+\ri I)v)|=\|Lv\|^2+\|v\|^2-2\re(\ri(Lv,v))\ge\|Lv\|^2+\|v\|^2
$$
because $\im(Lv,v)\ge0$.

If the operator $C$ defined by \rf{C} is a trace class operator, we say that $K$ is a {\it relatively trace class perturbation of $L$}.



It is well known and easy to verify that if $K$ is a relatively trace class perturbation of $L$, then it is a resolvent trace class perturbation of $L$, i.e. 
$$
(M+\ri I)^{-1}-(L+\ri I)^{-1}\in\bS_1.
$$

 The converse is not true which can easily be seen from Theorems \ref{x+iC}
and \ref{resLipC}.

\medskip

We are going to make three elementary (apparently well-known) observations. 

\medskip

\begin{lem}
\label{pervaya}
Suppose that 
\bay
\label{KcdA}
\|Kv\|\le c\|v\|+d\|Lv\|,\quad v\in{\rm Dom}(L),\quad\mbox{for some}\quad c>0\quad\mbox{and}\quad d\in(0,1)
\ey
(in this case we say that $K$ is {\it strictly dominated by $L$}). Then $K$ is also dominated by $L+tK$
for $t\in(0,1]$.
\end{lem}

\Pf
Indeed, let $v\in{\rm Dom}(L)$. We have
$$
\|Kv\|\le c\|v\|+d\|Lv\|=c\|v\|+d\|(L+tK)v-Kv\|\le c\|v\|+d\|(L+tK)v\|+dt\|Kv\|.
$$
It follows that
\bay
\label{cdK}
\|Kv\|\le(1-dt)^{-1}\big(c\|v\|+d\|(L+tK)v\|\big),
\ey
and so $K$ is dominated by $M$. $\bl$

\medskip

We proceed now to the case when $K$ is a {\it relatively compact perturbation of} $L$, i.e. the operator $K(A+\ri I)^{-1}$ is compact.

\begin{lem}
\label{vtoraya}
Suppose that $K$ is a {\it relatively compact perturbation of} $L$. Then for any positive number $d$, there exists a positive number $c$ such that
$$
\|Kv\|\le c\|v\|+d\|Lv\|,\quad v\in{\rm Dom}(L),
$$
and so, $K$ is strictly dominated by $M=L+K$.
\end{lem}

\Pf Let $\e>0$. By the assumption, the operator $K(L+\ri I)^{-1}$ is compact, and so it can be represented in the form $K(L+\ri I)^{-1}=T+R$, where $T$ is a finite rank operator and $\|R\|<\e$. 

Consider the operators $K_1$ and $K_2$ defined on ${\rm Dom}(A)$ by
$$
K_1=T(L+\ri I)\quad\mbox{and}\quad K_2=R(L+\ri I).
$$
Clearly,
$$
\|K_2v\|=\|R(L+\ri I)v\|\le\e\|(L+\ri I)v\|\le\e\|Lv\|+\e\|v\|,\quad v\in{\rm Dom}(L),
$$
and so, it suffices to show that for any positive number $d$, there exists a positive number $c$ such that
$$
\|K_1v\|\le c\|v\|+d\|Lv\|,\quad v\in{\rm Dom}(L).
$$
Since $T$ is a linear combination of rank one operators, it suffices to consider the case when 
$\rank T=\rank\big(K_1(L+\ri I)^{-1}\big)=1$. Suppose that
$$
K_1v=\big((L+\ri I)v,\f\big)\psi,\quad v\in{\rm Dom}(L).
$$
since $L$ is densely defined, $\f$ admit a representation $\f=\f_1+\f_2$, where 
$\f_1\in{\rm Dom}(L)$ and $\|\f_2\|<\e\|\psi\|^{-1}$. We have
\begin{align*}
K_1v&=\big((L+\ri I)v,\f_1\big)\psi+\big((L+\ri I)v,\f_2)\psi\big)\\[.2cm]
&=\big(v,(L+\ri I)\f_1)\psi+\big((L+\ri I)v,\f_2\big)\psi,
\quad v\in{\rm Dom}(L),
\end{align*}
and so
\begin{align*}
\|Kv\|&\le\|(L+\ri I)\f_1\|\cdot\|\psi\|\cdot\|v\|+\|\f_2\|\cdot\|\psi\|\cdot\|v\|+\|\f_2\|\cdot\|\psi\|\cdot\|Lv\|\\[.2cm]
&\le\big(\|(L+\ri I)\f_1\|\cdot\|\psi\|+\|\f_2\|\cdot\|\psi\|\big)\|v\|+\e\|Lv\|,
\quad v\in{\rm Dom}(L).\quad\bl
\end{align*}

\begin{lem}
\label{tret'ya}
Suppose that $K$ is a relatively trace class perturbation of  $L$.
Then $K$ is a relatively trace class perturbation of $M=L+K$.
\end{lem}

\Pf It is easy to see that
$$
K(L+K+\ri I)^{-1}-K(L+\ri I)^{-1}=-K(L+K+\ri I)^{-1}K(L+\ri I)^{-1}\in\bS_1.
$$
Indeed, by the assumption, $K(L+\ri I)^{-1}\in\bS_1$. On the other hand, by Lemmata \ref{pervaya} and
\ref{vtoraya}, the operator $K(L+K+\ri I)^{-1}$ is bounded. $\bl$

\medskip

\begin{lem}
\label{12max}
Suppose that $L$ is a maximal dissipative operator and $M$ is a dissipative operator such that 
$\Dom(M)\subset\Dom(L)$ and $K=M-L$ satisfies {\em\rf{KcdA}} with $d\in(0,\frac12)$. Then $M$ is a maximal dissipative operator.
\end{lem}

\Pf Let $\vk>0$. It is easy to see that
$$
\|(L+\ri\vk I )^{-1}\|\le\vk^{-1}.
$$
We have
$$
M+\ri\vk I =K+L+\ri\vk I=\big(I+K(L+\ri\vk I )^{-1}\big)(L+\ri\vk I ).
$$
Let us show that $(M+\ri\vk I)^{-1}$ is a bounded operator for some $\vk>0$. It suffices to show that $\|K(L+\ri\vk I )^{-1}\|<1$.
Let $v$ be a vector in our Hilbert space. By \rf{KcdA}, we have
$$
\|K(L+\ri\vk I )^{-1}v\|\le c\|(L+\ri\vk I )^{-1}v\|+d\|L(L+\ri\vk I )^{-1}v\|
$$
Clearly, 
$$
L(L+\ri\vk I )^{-1}v=v-\ri\vk(L+\ri\vk I )^{-1}v
$$
Thus,
$$
\|K(L+\ri\vk I )^{-1}v\|\le c\vk^{-1}\|v\|+d\|v\|+d\|v\|=(2+c\vk^{-1})d\|v\|.
$$
It suffices to select $\vk>0$ so that $(2+c\vk^{-1})d<1$. $\bl$

\medskip

{\bf Remark 1.} We conjecture that the conclusion of Lemma \ref{12max} should hold under the weaker assumption $d<1$, i.e. under the assumption that $K$ is strictly dominated by $L$. However, even if it could happen that the dissipative operator $M=L+K$ is not maximal, we can consider instead of $M$ its canonical maximal extension (see \S\:\ref{maks}). In this paper we consider one-parametric families $L_t=L+tK$, $0\le t\le1$. If $L_t$ is not maximal, we can consider instead its canonical maximal extension. We are going to keep the same notation $L_t$ for maximal the extension. 

\begin{cor}
Suppose that the operator $K$ is a relatively compact perturbation of $L$. Then the conclusion of the lemma holds.
\end{cor}

\medskip

{\bf Definition.}
A function $f$ in $\CAb$ is called an {\it analytic relatively operator Lipschitz function} if there is a positive number $k$ such that
\bay
\label{rOL}
\|f(M)-f(L)\|\le k\|(M-L)(L+\ri I)^{-1}\|
\ey
whenever $L$ and $M$ are are maximal dissipative operator operators such that $M-L$ is a relatively bounded perturbation of $L$.
{\it We denote by ${\rm ROL}_{\rm A}$ the class of such analytic relatively operator Lipschitz function on} $\R$.

\medskip

Note that in this definition it is natural to require that $f$ should be bounded. Indeed, if $f$ is a function defined  
on $\clos\C_+$ and we apply \rf{rOL} to dissipative operators on the one-dimensional space $\C$, we obtain
$$
|f(\mu)-f(\l)|\le\const|(\mu-\l)(\l+\ri)^{-1}|,\quad\l,~\mu\in\clos\C_+,
$$
which implies that $f$ should be bounded.

\medskip

{\bf Remark 2.} Obviously, if $f\in{\rm ROL}_{\rm A}$,  then $f$ is operator Lipschitz, and so, 
$f$ is differentiable on $\R$. Indeed, \rf{rOL} implies that
$$
\|f(M)-f(L)\|\le k\|(M-L)(L+\ri I)^{-1}\|\le k\|M-L\|\cdot\|(L+\ri I)^{-1}\|\le k\|M-L\|.
$$

\

\section{\bf Analytic Schur multipliers}
\setcounter{equation}{0}
\label{pred}

\

In this section we quote recent results of \cite{AP9} about analytic Schur multipliers. Recall that 
$\fM(\R^2)$ is the class of functions on $\R^2$ that are Schur multipliers with respect to arbitrary Borel spectral 
measures on $\R$. 

\medskip

{\bf Definition.} Let $f$ be a function that is bounded and analytic in $\C_+\times\C_+$ and suppose that
the limit
$$
\lim_{r\to0^+}f(x_1+\ri r,x_2+\ri r)\df f(x_1,x_2)
$$
exists for all $(x_1,x_2)\in\R^2$. We say that $f$ is an {\it analytic Schur multiplier on} $\C_+\times\C_+$ if
the boundary-value function $f$ on $\R^2$ belongs to $\fM(\R^2)$. We denote the class of analytic Schur multipliers by $\fM_{\rm A}$.

We need the following two results from \cite{AP9} on analytic Schur multipliers.

\medskip

{\bf I.} Put
$$
{\rm C}_{{\rm A},\be}(\C_+^2)=\big\{f\in H^\be(\C_+^2):~f(z,\cdot),~\:f(\cdot,z)\in{\rm C}_{{\rm A},\be}(\C_+)\quad\mbox{for all}\quad z\in\clos\C_+\big\}.
$$
Suppose that $\Phi\in{\rm C}_{{\rm A},\be}(\C_+^2)$ and $\Phi\in\fM_{\rm A}$. Then $f$ belongs to the Haagerup tensor product ${\rm C}_{{\rm A},\be}(\C_+)\otimes_{\rm h}{\rm C}_{{\rm A},\be}(\C_+)$. Moreover,
$$
\|\Phi\|_{{\rm C}_{{\rm A},\be}(\C_+)\otimes_{\rm h}{\rm C}_{{\rm A},\be}(\C_+)}
=\big\|\Phi \big| \R^2\big\|_{\fM(\R^2)}.
$$

\medskip

{\bf II.} Let
\bay
\label{CZC+}
\CAr\df\big\{f\in\CAb:~\:\mbox{the limit}~\:\lim_{|\z|\to\be}f(\z)~\:\mbox{exists}\big\} 
\ey
and
$$
{\rm C}_{{\rm A}}(\C_+^2)=\big\{f\in H^\be(\C_+^2):~f(z,\cdot), ~\:f(\cdot,z)\in{\rm C}_{{\rm A}}(\C_+)\quad\mbox{for all}\quad z\in\clos\C_+\big\}.
$$
Suppose that $\Phi\in{\rm C}_{{\rm A}}(\C_+^2)$ and $\Phi\in\fM_{\rm A}$. Then $f$ belongs to the Haagerup tensor product ${\rm C}_{{\rm A},\be}(\C_+)\otimes_{\rm h}{\rm C}_{{\rm A},\be}(\C_+)$. Moreover,
$$
\|\Phi\|_{{\rm C}_{{\rm A},\be}(\C_+)\otimes_{\rm h}{\rm C}_{{\rm A},\be}(\C_+)}
=\|\Phi \big| \R^2\|_{\fM(\R^2)}.
$$

\

\section{\bf A representation of operator differences\\ in terms of double operator integrals}
\setcounter{equation}{0}
\label{pred}

\

In this section we obtain a representation of the operator differences $f(M)-f(L)$ in terms of double operator integrals in the case when the function $\dg_\flat f$ on $\R^2$ defined by
\bay
\label{D(y+i)}
(\dg_\flat f)(x,y)\df
\left\{
\begin{array}{ll}\displaystyle\frac{f(x)-f(y)}{x-y}(y+\ri),&x\ne y\\[.4cm]
f'(x)(x+\ri),&x=y,
\end{array}
\right.
\ey
is a Schur multiplier with respect to arbitrary spectral (or semi-spectral) measures. Recall that the latter means that the function $\dg_\flat f$ belongs to the Haagerup tensor product 
$\CAb\otimes_{\rm h}\CAb$.

Note that similar results for functions of self-adjoint operators were obtained in \cite{AP6} and \cite{AP7}.

\begin{thm}
\label{dgbemol}
Let $L$ and $M$ be maximal dissipative operators with $M-L$ being a relatively bounded perturbation of $L$. Suppose that $f$ is a Lipschitz function in $\clos\C_+$,  analytic in $\C_+$ and such that
$$
\dg_\flat f\in  \CAb\otimes_{\rm h}\CAb.
$$
Then the following formula holds
\bay
\label{opra}
f(M)-f(L)=\iint_{\R\times\R}\frac{f(x)-f(y)}{x-y}(y+\ri)\,{\rm d}\E_{M}(x)(M-L)(L+\ri I)^{-1}\,{\rm d}\E_L(y),
\ey
and so $f$ must be relatively operator Lipschitz.
\end{thm}

\Pf It follows from Lemma 8.3 of \cite{AP7} and from {\bf II} of \S\:\ref{pred} that there are sequences $\{\f_n\}_{n\ge0}$ and 
$\{\psi_n\}_{n\ge0}$ of  functions in the space $\CAr$ (see \rf{CZC+}) 
such that the following statements hold:

(a)
$\sum\limits_{n\ge0}|\f_n|^2\le\const$ everywhere on $\C_+$;

(b)
$\sum\limits_{n\ge0}|\psi_n|^2\le\const$ everywhere on $\C_+$;

(c)
\bay
\label{Dflatfnpsin}
(w+\ri)(\dg f)(z,w) =\sum\limits_{n\ge0}\f_n(w)\psi_n(z)\quad z,~w\in\clos\C_+,
\ey
 where the divided difference $\dg f$ is defined by
 $$
(\dg f)(z,w)\df
\left\{
\begin{array}{ll}\displaystyle\frac{f(z)-f(w)}{z-w},&z\ne w\\[.4cm]
f'(w),&z=w.
\end{array}
\right.
 $$
 We have
 \begin{align}
 \label{ogrQ}
 \iint_{\R\times\R}(\dg_\flat f)(x,y)\,{\rm d}E_{M}(x)(M&-L)(L+\ri I)^{-1}\,{\rm d}E_L(y)\nonumber\\[.2cm]
 &=\sum_{n\ge0}\f_n(M)(M-L)(L+\ri I)^{-1}\psi_n(L)\df Q,
 \end{align}
 where $Q$ is a bounded linear operator.
 
%

 
Clearly,
\bay
\label{limity}
 \iint\limits_{\R\times\R}\left( f(y)(x+\ri)^{-1}-f(x)(x+\ri)^{-1}\right)\,\rd\E_M(y)\rd\E_L(x)
 =\big(f(M)-f(L)\big)(L+\ri I)^{-1}.
\ey
 On the other hand, for $z,~w\in\clos\C_+$, we have
 \begin{multline}
 \label{Oz,w}
 f(w)(z+\ri)^{-1}-f(z)(z+\ri)^{-1}=\sum_n\f_n(w)(w-z)(w+\ri)^{-1}(z+\ri)^{-1}\psi_n(z)\\[.2cm]
 =\sum_n\f_n(w)w(w+\ri)^{-1}(z+\ri)^{-1}\psi_n(z)\\[.2cm]
 -\sum_n\f_n(w)(w+\ri)^{-1}z(z+\ri)^{-1}\psi_n(z)
 \df\O(z,w).
 \end{multline}
Thus,
 \begin{multline*}
 \iint_{\R\times\R}\sum_n\f_n(y)y(y+\ri)^{-1}(x+\ri)^{-1}\psi_n(x)\,\rd\E_M(y)\,\rd\E_L(x)\\[.2cm]
 =\sum_n\f_n(M)\big(M(M+\ri I)^{-1}\big)(L+\ri I)^{-1}\psi_n(L)\\[.2cm]
 =\sum_n\f_n(M)\big(I-\ri(M+\ri I)^{-1}\big)(L+\ri I)^{-1}\psi_n(L).
 \end{multline*}
 Similarly,
 \begin{multline*}
 \iint_{\R\times\R}\left(
 \sum_n\f_n(y)(y+\ri)^{-1}x(x+\ri)^{-1}\psi_n(x)\right)\,\rd\E_M(y)\rd\E_L(x)\\[.2cm]
 =\sum_n\f_n(M)(M+\ri I)^{-1}(I-\ri(L+\ri I)^{-1})\psi_n(L).
 \end{multline*}
 It follows that
 \begin{multline*}
 \iint_{\R\times\R}\O(x,y)\,\rd\E_M(y)\rd\E_L(x)=
 \sum_n\f_n(M)(I-(M+\ri I)^{-1})(L+\ri I)^{-1}\psi_n(L)\\[.2cm]
 -\sum_n\f_n(M)(M+\ri I)^{-1}(I-(L+\ri I)^{-1})\psi_n(L)\\[.2cm]
 =\sum_n\f_n(M)\big((L+\ri I)^{-1}-(M+\ri I)^{-1}\big)\psi_n(L)\\[.2cm]
 =\sum_n\f_n(M)(M+\ri I)^{-1}(M-L)(L+\ri I)^{-1}\psi_n(L),
 \end{multline*}
 where the function $\O$ is defined in \rf{Oz,w}. Combining this equality with \rf{Oz,w} and \rf{limity}, we find that
$$
\big(f(M)-f(L)\big)(L+\ri I)^{-1}=\sum_n\f_n(M)(M+\ri I)^{-1}(M-L)(L+\ri I)^{-1}\psi_n(L).
$$
Multiplying this identity on the right by $(L+\ri I)$, we obtain
$$
f(M)-f(L)=\sum_n\f_n(M)(M+\ri I)^{-1}(M-L)\psi_n(L)=Q
$$
on $\Range(L+\ri I)=\Range L$ and since the operators on the left and on the right are bounded, it follows that the equality holds on the whole Hilbert space.  $\bl$
 
\medskip
 
\begin{cor}
Suppose that $f$ satisfies the hypotheses of Theorem {\em\ref{dgbemol}}. Let $L$ and $M$ be maximal dissipative operators such that $M-L$ is a relatively trace class perturbation of $L$. Then
$f(M)-f(L)\in\bS_1$.
\end{cor}

\

The converse of Theorem \ref{dgbemol} is also true. Indeed, the following result holds.

\begin{thm}
\label{veriobr}
Let $f\in\ROL_{\rm A}$. Then $\dg_\flat f\in\CAr\otimes_{\rm h}\CAr$, and so 
{\em\rf{opra}} holds whenever $L$ and $M$ are maximal dissipative operators such that $M-L$ is a relatively bounded perturbation of $L$.
\end{thm}

\Pf By Theorem 6.1 of \cite{AP7}, $\dg_\flat f\in\fM(\R^2)$, and so by {\bf II} of \S\:\ref{pred}, the function
$\dg_\flat f$ belongs to $\CAr\otimes_{\rm h}\CAr$ which implies the result. $\bl$

\

\section{\bf Various characterizations of analytic relatively operator Lipschitz functions}
\setcounter{equation}{0}
\label{khara}

Recall that $\OL_{\rm A}$ denote the space of operator Lipschitz functions on $\clos\C_+$ analytic in $\C_+$.
Note that $\OL_{\rm A}=\COL(\C_+)$, the space of commutator Lipschitz functions on $\C_+$, see \cite{APol}.

We say that a function $f$ on $\clos\C_+$ is a multiplier of the space $\OL_{\rm A}$  if
\bay
\label{mola}
\f\in\OL_{\rm A}\quad\Longrightarrow\quad \f f\in\OL_{\rm A}.
\ey
We denote by $\MOLA$ the class of all multipliers of $\OL_{\rm A}$.
By the norm $\|\f\|_{\MOLA}$ of a function $\f$ in $\MOLA$ we mean the norm of the multiplication operator 
$\f\mapsto\f f$ on the space $\OL_{\rm A}$.

\begin{thm}
\label{harotopL}
\label{x+iC}
Let  $f$ be a continuous function on $\clos\C_+$ that is analytic in $\C_+$ and differentiable on $\R$.
The following are equivalent:

{\em(a)} $f$ is an analytic relatively operator Lipschitz function;

{\em(b)} 
$(z+\ri)({\dg f})(z,w)\in\fM_{\rm A}$; 

{\em(c)}
$(z+\ri)f(z)\in\OLA$;

{\em(d)} $f\in\MOLA$.
\end{thm}

\Pf 
The implication (b)$\Rightarrow$(a) follows from Theorem \ref{dgbemol}.
The implication (a)$\Rightarrow$(b) follows from Theorem \ref{veriobr}.

It remains to prove that (b)$\Leftrightarrow$(c)$\Leftrightarrow$(d).

Clearly, (b) implies the condition 

$({\rm b}')$ $(x+\ri)({\dg f})(x,y)\in\fM(\R^2)$,

and $(c)$  implies the condition

$({\rm c}')$ $(x+\ri)f(x)\in\OL(\R)$.

 The implication $({\rm b}')\Rightarrow$(b)
 follows from the results of \cite{AP9}. 
 The implication $({\rm c}')\Rightarrow$(c) is a consequence of Theorem 3.9.5 in \cite{APol}.
 
  The equivalence (b)$\Leftrightarrow$(c) follows now from the equivalence $({\rm b}')\Leftrightarrow({\rm c}')$,
 which was proved in  Theorem 4.1  of \cite{AP7}.

The implication (d)$\Rightarrow$(c) is trivial.

It remains to prove that  (c)$\Rightarrow$(d). Let $(z+\ri)f(z)\in\OLA(\C_+)$.
Then $(x+\ri)f(x)\in\OL(\R)$, and $f\in\MOL$ by Theorem  4.1 in \cite{AP7}.
Thus, 
$$
\f\in\OL(\R)\quad\Longrightarrow\quad \f f\in\OL(\R).
$$
To complete the proof, we observe that $f$ is continuous on $\clos\C_+$ and analytic on $\C_+$. $\bl$

\begin{cor}
Let $f$ be a function of class  $\CAb$. Suppose that the function
$$
x\mapsto (x+\ri)f(x),\quad x\in\R,
$$
belongs to the Besov class $B^1_{\be,1}(\R)$. Then $f$ is an analytic relatively operator Lipschitz function.
\end{cor}

\Pf This is an immediate consequence of Theorem \ref{harotopL} and the fact that the Besov class 
$B^1_{\be,1}(\R)$ consists of operator Lipschitz functions, see \cite{Pe1} and \cite{Pe3}. $\bl$

\medskip

Let us proceed now to various descriptions of the class of analytic resolvent Lipschitz functions.

\begin{thm}
\label{resLipC} 
Let  $f$ be a continuous function on $\clos\C_+$ that is analytic in $\C_+$ and differentiable on $\R$.
The following are equivalent:

{\em(a)} $f\in\OLAr$, i.e. $f$ is an an analytic resolvent Lipschitz function;

{\em(b)}
$(z+\ri)(w+\ri)(\dg f)(z,w)\in\frak M_{\rm A}(\C_+^2)$;

{\em(c)}
$(z+\ri)^2(f(z)-c)\in\OLA(\C_+)$ for some $c\in\C$;

{\em(d)}
$(z+\ri)(f(z)-c)\in\MOLA$ for some $c\in\C$.

{\em(e)} the function $\f$ on $\clos\dd\setminus\{1\}$ defined by 
\bay
\label{peresazhennaya}
\f(\z)\df f\left(\ri\frac{1+\z}{1-\z}\right),\quad\z\in\clos\dd\setminus\{1\},
\ey
extends to a function of class $\OL_{\rm A}(\dd)$
of operator Lipschitz function on $\clos\dd$ analytic in $\dd$.

If $c\in\C$ satisfies {\em(b)} or {\em(c)}, then $c$ must be equal to $\lim_{|z|\to\be}f(z)$
\end{thm}

\Pf \Pf First we prove that (a)$\Rightarrow$(e). By (a), we have
\bay
\label{lipAB}
\|f(M)-f(L)\|\le\const\|(M+\ri I)^{-1}-(L+\ri I)^{-1}\|
\ey
for all maximal dissipative operators $L$ and $M$.

Considering dissipative operators on the one-dimensional space, we find that
\rf{lipAB} implies that
$$
|f(w)-f(z)|\le\const|(w+\ri)^{-1}-(z+\ri)^{-1}|,\quad z,~w\in\clos\C_+,
$$
and so, the limit $\lim_{|z|\to\be}f(z)$ exists. Thus, the function $\f$ defined by \rf{peresazhennaya}
extends to a function in the disc-algebra.

Let $U$ and $V$ be the Cayley transforms of $L$ and $M$, i.e.
\bay
\label{UUAA}
U\df(L-\ri I)(L+\ri I)^{-1}\quad \text{and} \quad V\df(M-\ri I)(M+\ri I)^{-1}.
\ey
Clearly, $U$ and $V$ are unitary operators if $L$ and $M$ are self-adjoint operators.
Now applying \rf{lipAB} to the bounded self-adjoint operators  $L$ and $M$,
we find that
\bay
\label{lipUV}
\|\f(V)-\f(U)\|\le\const\|V-U\|
\ey
for all unitary operators $U$ and $V$ with $1\notin\s(U)$
and $1\notin\s(V)$.

To obtain \rf{lipUV} for arbitrary unitary operators $U$ and $V$, we can construct two sequences
of unitary operators $\{U_n\}_{n\ge1}$ and $\{V_n\}_{n\ge1}$ such that
$1\notin\s(U_n)$ and $1\notin\s(V_n)$ for all $n\ge1$,
$U_nU=UU_n$ and $V_nV=VV_n$ for all $n\ge1$,
$\lim_{n\to\be}\|U_n-U\|=0$ and
$\lim_{n\to\be}\|V_n-V\|=0$.

It remains to observe that the function $\f$ extends to a function
in the disk-algebra.

Let us prove that  (e)$\Rightarrow$(a). By Theorem 3.9.9 of \cite{APol}, inequality \rf{lipUV} holds
for arbitrary contractions $U$ and $V$. We want to prove inequality \rf{lipAB}
for the maximal dissipative operators $L$ and $M$. It remains to observe that inequality \rf{lipAB}
is equivalent to inequality \rf{lipUV} for the contractions $U=(L-\ri I)(L+\ri I)^{-1}$ and
$V=(M-\ri I)(M+\ri I)^{-1}$.

The equivalence of (e) and (c) follows from Theorem 3.10.2 in \cite{APol}.

In what follows we need the following identity:
\begin{align}
\label{iziw}
\frac{(z+\ri)^2f(z)-(w+\ri)^2f(w)}{z-w}&=(z+\ri)(w+\ri)(\dg f)(z,w)\nonumber\\[.2cm]
&+(z+\ri)f(z)+(w+\ri)f(w).
\end{align}
 
Let us show that (b)$\Rightarrow$(c).
Suppose that $(z+\ri)(w+\ri)(\dg f)(z,w)\in\frak M_{\rm A}(\C_+^2)$. Then \lb$(z+\ri)(\dg f)(z,w)\in\frak M_{\rm A}(\C_+^2)$.
It follows from the previous theorem that $(z+\ri)f(z)\in\OLA(\C_+)$. Thus,
the limit $c\df\lim_{|z|\to\be}f(z)$ exists. Without loss of generality we may assume that $c=0$.
To deduce from \rf{iziw} the inclusion $(z+\ri)^2f(z)\in\OLA(\C_+)$, it suffices to observe that
the function $(z+\ri)f(z)$ is bounded. This follows from the fact that $\lim_{|z|\to+\infty}f(x)=0$ and from the fact 
that the function $(z+\ri)^2(\dg f)(z,z)=(z+\ri)^2f'(z)$ is bounded.

Let us prove that (c)$\Rightarrow$(b).
Suppose that $(z+\ri)^2(f-c)\in\OLA(\R)$ for some $c\in\C$. Let us show that $(z+\ri)(w+\ri)(\dg f)(z,w)\in\frak M(\C_+^2)$. Again, we may assume that $c=0$. 
The result follows from the fact that the function $(z+\ri)f(z)$ is bounded and from identity \rf{iziw}. 

The equivalence of (c) and (d) is a consequence of in the statement of Theorem \ref{x+iC}.

Finally, it is obvious that (c) can hold only for one value of $c$. The same is true for (d).
$\bl$

\

\section{\bf Differentiaton in the parameter}
\setcounter{equation}{0}
\label{differ}

\

As before, $L$ is a maximal dissipative operator and $M$ is a dissipative operator such that $K=M-L$ is a relatively bounded perturbation of $L$. In this section for a function for a relatively operator function $f$ analytic in $\C_+$ we consider the problem of the differentiability of the map $t\mapsto f(L+tK)$. Obviously, for sufficiently small positive $t$, the operator $tK$ is strictly dominated by $L$.

\begin{thm}
\label{sil'no}
Let $L$ be a maximal dissipative operator and let $M$ be a dissipative operator such that $K=M-L$
is a relatively bounded perturbation of $L$.
Suppose that $f\in\ROL_{\rm A}$. Then the function
$t\mapsto f(L_t)\df f(L+tK)$ is differentiable at 0 in the strong operator topology and
\bay
\label{proizvodnaya}
\frac{\rm d}{{\rm d}t}f(L_t)\Big|_{t=0}=
\iint_{\R\times\R}\frac{f(x)-f(y)}{x-y}(y+\ri)\,{\rm d}\E_{L}(x)C\,{\rm d}\E_{L}(y)
\ey
where the operator $C$ and is defined by {\em\rf{C}}.
\end{thm}

\Pf By Theorem \ref{veriobr} the function $\dg_\flat f$ belongs to the Haagerup tensor product 
$\CAr\otimes_{\rm h}\CAr$, and so there are sequences $\f_n$ and $\psi_n$ in $\CAr$ such that
$$
\sup_{z\in\clos\C_+}\sum_n|\f_n(z)|^2<\be,\quad\sup_{z\in\clos\C_+}\sum_n|\psi_n(z)|^2<\be
$$
and
$$
\dg_\flat f(z,w)=\sum_n\f_n(z)\psi_(w),\quad z,~w\in\clos\C_+.
$$

Let $K=M-L$. Consider the parametric family $L_t=L+K_t$, where $t$ is a sufficiently small positive number. Then $L_t$ is a dissipative operator. Denote by $\E_t$ the semi-spectral measure of $L_t$. We have
\begin{align*}
\frac1t(f(L_t)-f(L))&=
\iint_{\R\times\R}\frac{f(x)-f(y)}{x-y}(y+\ri)\,{\rm d}\E_t(x)K(L+\ri I)^{-1}\,{\rm d}\E_0(y)\\[.2cm]
&=\sum_{n\ge0}\f_n(L_t)K(L+\ri I)^{-1}\psi_n(L),
\end{align*}
We have
\bay
\label{neprer-taki}
\lim_{t\to\0}\|\f_n(L_t)-\f_n(L)\|=0.
\ey
This follows from the fact that function $\f$ in the class 
$\CAr$ is uniformly operator continuous, i.e.
for each $\e>0$, there exists $\d>0$ such that 
\bay
\label{anopco}
\|\f_n(L)-\f_n(M)\|<\e
\ey
whenever $M$ is a maximal dissipative operator satisfying $\|L-M\|<\d$.
The latter is a consequence of Theorem 7.2 of \cite{AP1} and an analog of Theorem 8.1 of \cite{APH} for maximal dissipative operators.

We have to prove that
$$
\lim_{t\to0}\sum_{n\ge0}\f_n(L_t)K(L+\ri I)^{-1}\psi_n(L)=\sum_{n\ge0}\f_n(L)K(L+\ri I)^{-1}\psi_n(L)
$$
in the strong operator topology. Thus, we need to prove that for an arbitrary vector $u$,
$$
\lim_{t\to0}\sum_{n\ge0}(\f_n(L_t)-\f_n(L))K(L+\ri I)^{-1}\psi_n(L)u=\bf 0
$$
where the series is understood in the sense of  the weak topology of $\mathscr H$ while the limit is taken
in the norm of $\mathscr H$. Assume that $\|u\|=1$.
Let $u_n\df K(L+\ri I)^{-1}\psi_n(L)u$. We have
\begin{align*}
\sum_{n\ge0}\|u_n\|^2&\le\|K(L+\ri I)^{-1}\|^2\sum_{n\ge0}\|(\psi_n(L)u\|^2
\le\|K(L+\ri I)^{-1}\|^2\sum_{n\ge0}\|(\psi_n(A)u\|^2
\\[.2cm]
&=\|K(L+\ri I)^{-1}\|^2\sum_{n\ge0}(|\psi_n|^2(A)u,u)
\le\sup_{z\in\clos\C_+}\left(\sum_n|\f_n(z)|^2\right)\|K(L+\ri I)^{-1}\|^2
\end{align*}
where $A$ is a resolvent self-adjoint dilation of $L$.

Let $\e>0$. We can choose a natural number $N$ such that $\sum_{n>N}\|u_n\|^2<\e^2$.
 Then it
follows from Lemma 3.5.9 in \cite{APol} that
\begin{multline}
\label{khvost}
\Big\|\sum_{n>N}(\f_n(L_t)-\f_n(L))u_n \Big\|\le\Big\|\sum_{n>N}\f_n(L_t)u_n\Big\|\nonumber\\[.2cm]
+\Big\|\sum_{n>N}\f_n(L)u_n\Big\|
<2\e\sup_{z\in\clos\C_+}\left(\sum_n|\f_n(z)|^2\right)^{1/2}
\end{multline}
for all $t\in[0,1]$. By \rf{neprer-taki},
$$
\left\|\sum_{n=0}^N(\f_n(A+tK)-\f_n(A))u_n \right\|<\e
$$
for all $t$ sufficiently close to zero. It follows that
$$
\Big\|\sum_{n\ge0}(\f_n(L_t)-\f_n(L))u_n\Big\|<
\left(2\sup_{z\in\clos\C_+}\left(\sum_n|\f_n(z)|^2\right)^{1/2}+1\right)\e
$$
for all $t$ sufficiently close to zero. $\bl$

\begin{thm}
Let $L$ be a maximal dissipative operator and let $M$ be a dissipative operator such that $K=M-L$ is 
strictly dominated by $L$. Suppose that $f\in\ROL_{\rm A}$. Then for $t\in(0,1)$, the following identity holds
$$
\frac{\rm d}{{\rm d}s}f(L_s)\Big|_{s=t}=
\iint_{\R\times\R}\frac{f(x)-f(y)}{x-y}(y+\ri)\,{\rm d}\E_{L_t}(x)K(L_t+\ri I)^{-1}\,{\rm d}E_{L_t}(y),
$$
in the strong operator topology, where $L_t\df L+tK$.
\end{thm}

\Pf The result is an immediate consequence of Theorem \ref{sil'no} and Lemma \ref{pervaya}. $\bl$

\

\section{\bf The trace formula in the case of relatively trace class perturbations}
\setcounter{equation}{0}
\label{SlROLA}

\

In this section we obtain a trace formula for the trace of $f(M)-f(L)$, where 
$L$ be a maximal dissipative operator and let $M$ be a dissipative operator such that $M-L$ is a relatively trace class perturbation of $L$. We prove that the pair
$\{L,M\}$ has spectral shift function $\bs\xi$ such that 
$$
\int_\R\frac{|\bs\xi(t)|}{1+|t|}\,\rd t<\be.
$$
We also show in this section that the maximal class of functions $f$, for which the trace formula holds with arbitrary pairs $\{L,M\}$ of maximal dissipative operators with $M-L$ being a relatively trace class perturbation of $L$ coincides with $\ROL_{\rm A}$.

Put $L_t\df L+t(M-L)$, $t>0$.
It follows from Lemmata \ref{vtoraya}, \ref{tret'ya} and \ref{12max} that $M-L$ is a relatively trace class perturbation of $L_t$ and $L_t$ is a maximal dissipative operator.

\begin{thm}
\label{flasle}
Let $L$ be a maximal dissipative operator and let $M$ be a dissipative operator such that $M-L$ is a relatively trace class perturbation of $L$. Then there exists a measurable function $\bs\xi$ such that
\bay
\label{vpervoi}
\int_\R\frac{|\bs\xi(t)|}{1+|t|}\,\rd t<\be
\ey
and the trace formula 
\bay
\label{nasledie}
\trace\big(f(M)-f(L)\big)=\int_\R f'(t)\bs\xi(t)\,\rd t
\ey
holds for an arbitrary function $f$ in $\ROL_{\rm A}$. Moreover, $\ROL_{\rm A}$ is the maximal class of functions $f$, for which trace formula {\em\rf{nasledie}} for all such pairs $\{L,M\}$.
\end{thm}

A function $\bs\xi$ on $\R$ is called a {\it spectral shift function for the pair} $\{L,M\}$ if it satisfies \rf{vpervoi} and
\rf{nasledie}.

We need the following lemma, which is similar to Lemma 9.2 in \cite{AP7}.

\begin{lem}
\label{neprogr}
The function $t\mapsto K(L_t+\ri I)^{-1}$ is continuous on $[0,1]$ in the norm of $\bS_1$, and so
$\|K(L_t+\ri I)^{-1}\|_{\bS_1}\le\const$, $0\le t\le1$. 
\end{lem}

\Pf We have
$$
K(L_t+\ri I)^{-1}=K(L+\ri I)^{-1}(L+\ri I)(L_t+\ri I)^{-1}.
$$
Clearly, it suffices to show that the function 
$$
t\mapsto\big((L+\ri I)(L_t+\ri I)^{-1}\big)^{-1}=(L_t+\ri I)(L+\ri I)^{-1}
$$ 
is continuous. This follows immediately from the following obvious equality
$$
(L_t+\ri I)(L+\ri I)^{-1}=I+tK(L+\ri I)^{-1}.\quad\bl
$$

\medskip

\Pf Consider the one-parametric family of operators $L_t\df L+tK$, $0\le t\le1$. 

By Theorem \ref{sil'no},
\bay
\label{LQt}
\!\!\frac{\rm d}{{\rm d}s}\big(f(L_s)-f(L_t)\big)\Big|_{s=t}=
\iint_{\R\times\R}\!\!\frac{f(x)-f(y)}{x-y}(y+\ri)\,{\rm d}\E_{L_t}(x)K(L_t+\ri I)^{-1}\,{\rm d}\E_{L_t}(y),
\ey
where the derivative on the left is understood in the strong operator topology.

Let $Q_t$ be the operator on the right-hand side of \rf{LQt}. Since the function  $\dg_\flat f$ is a Schur multiplier for arbitrary semi-spectral measures, it follows that
$Q_t\in\bS_1$ for every $t$ in $[0,1]$ and $\sup_t\|Q_t\|_{\bS_1}<\be$, see Lemma \ref{neprogr}.

Since the function $t\mapsto Q_t$ is continuous in the norm of $\bS_1$ by Lemma \ref{neprogr},
we obtain
$$
f(L+K)-f(A)=\int_0^1Q_t\,{\rm d}t
$$
in the sense of integration of continuous functions and
$$
\trace\big(f(A+K)-f(A)\big)=\int_0^1\trace Q_t\,{\rm d}t.
$$
We need the following lemma.

\begin{lem}
\label{vychsl}
Let $L$ be a maximal dissipative operator and let $R$ be a trace class operator. 
Suppose that $f\in\ROL_{\rm A}$.
Then
$$
\trace\iint_{\R\times\R}(\dg_\flat f)(x,y)\E_L(x)R\,{\rm d}\E_L(y)
=\int_\R f'(x)(x+\ri)\,\rd\nu(x),
$$
where $\nu$ is the complex Borel measure on $\R$ defined by $\nu(\D)=\trace(\E(\D)R)$ for a Borel subset $\D$ of $\R$.
\end{lem}

\Pf Let  $\{\f_n\}_{n\ge0}$ and 
$\{\psi_n\}_{n\ge0}$ be sequences of  functions in $\CAr$ that satisfy conditions a), b) and c) of Theorem \ref{dgbemol}. We have
\begin{multline*}
\trace\iint_{\R\times\R}(\dg_\flat f)(x,y)\E_L(x)R\,{\rm d}\E_L(y)=\sum_n\trace\f_n(L)R\f_n(L)
\\[.2cm]
=\trace\left(R\sum_n(\f_n\psi_n)(L)\right)=
\trace\left(R\int_\R\left(\sum_n(\f_n\psi_n)\right)(x)\,\rd\E_L(x)\right)\\[.2cm]
=\trace\left(R\int_\R(\dg_\flat f)(x,x)\,\rd\E_L(x)\right)=\trace\left(R\int_\R f'(x)(x+\ri)\,\rd\E_L(x)\right)\\[.2cm]
=\int_\R f'(x)(x+\ri)\,\rd\nu(x).\quad\bl
\end{multline*}

Let us return to the proof of Theorem \ref{flasle}.

We apply Lemma \ref{vychsl} to the dissipative operator $L_t$ and the trace class operator $K(L+\ri)^{-1}$
and find that
$$
\trace Q_t=\int_\R f'(x)(x+\ri)\,\rd\nu_t(x),
$$
where $\nu_t$ is the complex measure defined by $\nu_t(\D)=\trace\big(\E_t(\D)K(L_t+\ri I)^{-1}\big)$ for a Borel subset $\D$ of $\R$. Clearly $\|\nu_t\|\le\|K(L_t+\ri I)^{-1}\|_{\bS_1}$. 

We consider the space $\mM(\R)$ of complex Borel measures on $\R$ as the dual space to the space
${\rm C}_0(\R)$ of continuous functions on $\R$ vanishing at infinity with respect to the pairing
\bay
\label{sparim}
\langle h,\mu\rangle=\int_\R h\,\rd\mu.
\ey
Consider the subspace ${\rm C}_0(\R)\cap\CAr$ of ${\rm C}_0(\R)$ whose dual space can be identified naturally with the quotient space $\mM(\R)/\big({\rm C}_0(\R)\cap\CAr\big)^\perp$ with respect to the pairing
\ref{sparim}.

Given a measure $\mu$ in $\mM(\R)$, we denote by $\dot\mu$ the coset in 
$\mM(\R)/\big({\rm C}_0(\R)\cap\CAr\big)^\perp$, which corresponds to $\mu$.

Let us show that the function $t\mapsto\dot\nu_t$ is continuous in the weak-$*$ topology of 
\lb$\mM(\R)/\big({\rm C}_0(\R)\cap\CAr\big)^\perp$.
To this end we establish the following lemma. 

\begin{lem}
\label{Lnu}
Let $L$ be a maximal dissipative operator and let $R\in\bS_1$
Consider the complex Borel measure $\mu$ on $\R$ defined by
\bay
\label{borD}
\mu(\D)=\trace\big(\E_L(\D)R\big),\quad\D\subset\R.
\ey
Then
\bay
\label{anatozh}
\int_\R h\,\rd\mu=\trace \big(h(L)R\big)
\ey
for an arbitrary function $h$ in ${\rm C}_0(\R)\cap\CAr$.
\end{lem}

\Pf Consider a resolvent self-adjoint dilation $A$ of $L$ and observe that \rf{anatozh} is equivalent to the identity
\bay
\label{tozhdh}
\int_\R h\,\rd\mu=\trace \big(h(A)R\big).
\ey
In fact, \rf{tozhdh} holds for an arbitrary continuous function on $\R$ satisfying
$\lim_{|t|\to\be}f(t)=0$. Indeed, \rf{borD} means that \rf{tozhdh} holds characteristic functions. It suffices to approximate the function $h$ in \rf{tozhdh} by linear combinations of characteristic functions. $\bl$

\medskip

Let us return to the proof of Theorem \ref{flasle}. 
We want prove that the function $t\mapsto\dot\nu_t$ is  weak-$*$ continuous, i.e. the map
\bay
\label{hnut}
t\mapsto\int_\R h\,\rd\nu_t
\ey
is continuous for an arbitrary function $h$ in ${\rm C}_0(\R)\cap\CAr$.
Suppose first that $h\in{\rm C}_0(\R)\cap\CAr$. Then by Lemma \ref{Lnu},
$$
\int_\R h\,\rd\nu_t=\trace\big(h(L_t)K(L_t+\ri I)^{-1}\big).
$$
The continuity of the map \rf{hnut} follows now from Lemma \ref{neprogr} and from inequality \rf{anopco}.

Let $\nu$ be a measure in $\mM(\R)$ such that
$$
\dot\nu=\int_0^1\nu_t\,\rd t.
$$
We can also require that 
$$
\|\nu\|_{\mM(\R)}\le\max_{0\le t\le1}\|K(L_t+\ri I)^{-1}\|,
$$
which is finite by Lemma \ref{neprogr}.

Next, by Lemma \ref{neprogr}, the function $t\mapsto Q_t$ is continuous on $[0,1]$. 
On the other hand, the function $x\mapsto f'(x)(x+\ri)$ is a bounded function of the first Baire class (see the proof of Theorem 9.1 in \cite{AP7}). 

Thus, by Lemma 9.3 of  \cite{AP7},
\begin{align*}
\trace\big(f(M)-f(L)\big)&=\int_0^1\trace Q_t\,\rd t=\int_0^1\left(\int_\R f'(x)(x+\ri)\,\rd\nu_t(x)\right)\\[.2cm]
&=\int_\R f'(x)(x+\ri)\,\rd\nu(x).
\end{align*}

We have mentioned in \S\:\ref{otno} that the condition that $M-L$ is a relatively trace class perturbation of $L$ implies that
$$
(M+\ri I)^{-1}-(L+\ri I)^{-1}\in\bS_1. 
$$

It was shown in \S\:7 of \cite{MNP} that in this case the pair $\{L,M\}$ has a spectral shift function $\bs\eta$, such that
$$
\int_\R\frac{|\bs{\eta}(x)|}{1+x^2}\,{\rm d}x<\be
$$
and
$$
\trace\big(f(M)-f(L)\big)=\int_\R f'(x)\bs{\eta}(x)\,{\rm d}x
$$
whenever $f$ is a rational function with poles in $\C_-\df\{\z\in\C:~\im\z<0\}$. It follows that 
$$
\int_\R f'(x)\bs{\eta}(x)\,{\rm d}x=\int_\R f'(x)(x+\ri)\,{\rm d}\nu(x)
$$
for an arbitrary rational function $f$ with poles in $\C_-$. Thus,
$$
\int_\R\frac1{(\l-x)^{2}}\bs\eta(x)\,dx
=\int_\R\frac1{(\l-x)^{2}}\,d\nu(x),\quad\im\l<0.
$$

Consider the Radon complex measure $\mu$ defined by $d\mu=d\nu-\bs\eta\,d\m$, where $\m$ stands for Lebesgue measure on $\R$. Then
$$
\int_\R\frac1{(\l-x)^{2}}\,d\nu(x)=0,\quad\im\l<0.
$$
It can easily be deduced from
the brothers Riesz theorem (see e.g., Lemma 3.7 of \cite{MN} for details) that the complex Borel measure
 $\nu$ is absolutely continuous with respect to 
Lebesgue measure. We can define now the function $\bs{\xi}$ by
$$
\bs{\xi}(s)\df(s+\ri)\frac{{\rm d}\nu(x)}{{\rm d}x}(s).
$$
Clearly, $\bs\xi$ satisfies \rf{vpervoi} and equality \rf{nasledie} holds. $\bl$

\begin{thm}
\label{dismaksi}
Let $f$ be a function in $\CAr$ that is 
differentiable on $\R$ such that trace formula {\em\rf{nasledie}}
holds for arbitrary pairs $\{L,M\}$, where $L$ is a maximal dissipative operator and $M$ is a dissipative operator such that $M-L$ is a relatively trace class perturbation of $L$. Then $f\in\ROL_{\rm A}$.
\end{thm}

\Pf Since each self-adjoint operator is a maximal dissipative operator, the result follows from Theorem 9.4 in \cite{AP7}. $\bl$

\

\begin{footnotesize}

\
 
\noindent
\begin{tabular}{p{8cm}p{15cm}}
A.B. Aleksandrov & V.V. Peller \\
St.Petersburg State University & St.Petersburg State University \\
Universitetskaya nab., 7/9  & Universitetskaya nab., 7/9\\
199034 St.Petersburg, Russia & 199034 St.Petersburg, Russia \\
\\

St.Petersburg Department &St.Petersburg Department\\
Steklov Institute of Mathematics  &Steklov Institute of Mathematics  \\
Russian Academy of Sciences  & Russian Academy of Sciences \\
Fontanka 27, 191023 St.Petersburg &Fontanka 27, 191023 St.Petersburg\\
Russia&Russia\\
email: alex@pdmi.ras.ru&

email: peller@math.msu.edu

\end{tabular}
\end{footnotesize}

\end{document}

%% file: classstartscr.tex
\newcommand{\lb}{\linebreak}

\renewcommand{\d}{\delta}
\newcommand{\e}{\varepsilon}
\newcommand{\vk}{\varkappa}
\newcommand{\z}{\zeta}

\renewcommand{\l}{\lambda}

\newcommand{\s}{\sigma}

\newcommand{\f}{\varphi}
\renewcommand{\o}{\omega}

\newcommand{\D}{\Delta}

\renewcommand{\O}{\Omega}

\newcommand{\E}{{\mathscr E}}

\newcommand{\h}{{\mathscr H}}

\newcommand{\K}{{\mathscr K}}

\newcommand{\X}{{\mathscr X}}
\newcommand{\Y}{{\mathscr Y}}

\newcommand{\C}{{\Bbb C}}
\newcommand{\T}{{\Bbb T}}

\newcommand{\dd}{{\Bbb D}}
\newcommand{\R}{{\Bbb R}}
\newcommand{\Z}{{\Bbb Z}}

\newcommand{\0}{{\boldsymbol{0}}}

\newcommand{\bs}{\boldsymbol}

\newcommand{\m}{{\boldsymbol m}}
\newcommand{\bS}{{\boldsymbol S}}

\newcommand{\rf}[1]{(\ref{#1})}

\newcommand{\df}{\stackrel{\mathrm{def}}{=}}

\newcommand{\re}{\operatorname{Re}}
\newcommand{\spn}{\operatorname{span}}

\newcommand{\clos}{\operatorname{clos}}

\newcommand{\trace}{\operatorname{trace}}
\newcommand{\rank}{\operatorname{rank}}
\newcommand{\const}{\operatorname{const}}

\newcommand{\eeq}{\end{equation}}
\newcommand{\beq}{\begin{equation}}
\newcommand{\bay}{\begin{eqnarray}}
\newcommand{\ba}{\begin{align*}}
\newcommand{\ea}{\end{align*}}
\newcommand{\ey}{\end{eqnarray}}
\newcommand{\bey}{\begin{eqnarray*}}
\newcommand{\eey}{\end{eqnarray*}}

\newcommand{\be}{\infty}

\newcommand{\bl}{\blacksquare}

\newcommand{\Range}{\operatorname{Range}}

\newcommand{\Pf}{{\bf Proof. }}
\newcommand{\im}{\operatorname{Im}}
\renewcommand{\re}{\operatorname{Re}}

\newtheorem{thm}{\hspace{\parindent}Theorem}[section]

\newtheorem{cor}[thm]{\hspace{\parindent}Corollary}
\newtheorem{lem}[thm]{\hspace{\parindent}Lemma}